\input amstex
\documentstyle{amsppt}

\input label.def
\input degt.def
\def\Dg:{\endgraf{\bf Dg:}\enspace\ignorespaces}
\def\Kh:{\endgraf{\bf Kh:}\enspace\ignorespaces}

\newtheorem{Remark}\Remark\thm\endAmSdef

\def\dash{\item"\hfill--\hfill"}
\def\Dashes{\widestnumber\item{--}\roster}
\def\endDashes{\endroster}

{\catcode`\@11
\gdef\proclaimfont@{\sl}
\gdef\subsubheadfont@{\bf}
}

\def\paragraph{\subsubsection{}}
{\let\subhead\relax
\global\def\subsubhead#1\endsubsubhead{\subhead#1\endsubhead}}

\input epsf
\def\picture#1{\epsffile{#1-bb.eps}}

\def\ie{\emph{i.e.}}
\def\eg{\emph{e.g.}}
\def\cf{\emph{cf}}
\def\via{\emph{via}}
\def\etc{\emph{etc}}

\def\ZZ{\Z_2}
\def\MAX(#1,#2){B^0_{#1}(#2)}

\def\Cp#1{\Bbb P^{#1}}
\def\Rp#1{\Bbb P_\R^{#1}}
\def\tsmash#1{\smash{\tilde#1}}
\def\CR{C_\R}
\def\VR{V_\R}

\def\LR{L_\R}
\def\tL{\tsmash L}

\def\Lm{L_+}

\let\tq=q
\def\CpN{\Cp{N}}
\def\Cpr{\Cp{r}}
\def\RpN{\Rp{N}}
\def\Rpr{\Rp{r}}

\def\Sr{S^r}
\def\SS{\Omega}

\def\CZ{\Cal Z}
\def\CC{\Cal C}
\def\CK{\Cal K}
\def\CL{\Cal L}
\def\CO{\Cal O}

\def\tC{\tsmash C}

\def\imax{i_{\max}}
\def\Dmax{D_{\max}}
\def\ls|#1|{\mathopen|#1\mathclose|}
\def\oval{\frak o}
\def\arc{\frak a}
\def\oind#1{\ind #1}
\def\oind#1{i(#1)}
\def\omax#1{\operatorname{Hilb}(#1)}
\def\Hilb{\operatorname{Hilb}}
\def\ix#1{^{(#1)}}

\def\sset{\frak S}

\def\Spin{\operatorname{Spin}}
\def\pred{\operatorname{pred}}
\def\Arf{\operatorname{Arf}}
\def\Pic{\operatorname{Pic}}
\def\Card{\operatorname{Card}}
\def\corank{\operatorname{corank}}
\def\res{\operatorname{res}}
\def\ind{\operatorname{ind}}
\def\int{\operatorname{int}}

\def\dpth{\operatorname{dp}}

\topmatter

\author
Alex Degtyarev, Ilia Itenberg, Viatcheslav Kharlamov
\endauthor

\thanks
The second and the third authors are supported by
ANR-05-0053-01 grant of Agence Nationale de la Recherche (France)
and a grant of Universit\'{e} Louis Pasteur, Strasbourg.
\endthanks

\address
Bilkent University\endgraf\nobreak
06800 Ankara, Turkey
\endaddress

\email
degt\@fen.bilkent.edu.tr
\endemail

\address
Universit\'{e} Pierre et Marie Curie and
Institut Universitaire de France\endgraf\nobreak
Institut de Math\'{e}matiques de Jussieu\endgraf\nobreak
4 place Jussieu, 75005 Paris, France
\endaddress

\email
itenberg\@math.jussieu.fr
\endemail

\address
Institut de Recherche Math\'{e}matique Avanc\'{e}e,\endgraf\nobreak
Universit\'{e} de Strasbourg et CNRS,\endgraf\nobreak
7 rue Ren\'{e} Descartes 67084 Strasbourg Cedex, France
\endaddress

\email
kharlam\@math.u-strasbg.fr
\endemail

\title
On the number of
components of a complete intersection of real quadrics
\endtitle

\rightheadtext{components of intersection of real quadrics}

\abstract
Our main results concern complete intersections of three
real quadrics. We prove that the maximal number~$B^0_2(N)$ of
connected components that a regular complete intersection of three
real quadrics in $\Bbb{P}^N$ can have differs at most by one from
the maximal number of ovals of the submaximal depth $[(N-1)/2]$ of
a real plane projective
curve of degree $d=N+1$. As a consequence, we obtain
a lower bound \smash{$\frac14 N^2+O(N)$} and an upper bound
\smash{$\frac38 N^2+O(N)$} for $B^0_2(N)$.
\endabstract



\dedicatory
To Oleg Viro
\medskip
\rightline{\vbox{\hsize 80mm\eightpoint \leftskip0pt\raggedright\noindent
\bf viros / viro :\enspace\it
Terme gaulois d\'{e}signant ce qui est juste, vrai, sinc\`{e}re\dots
\smallskip\noindent\rightskip0pt \leftskip10mm plus2em\rm
X. Delamarre, Dictionnaire de la langue gauloise, Errance, Paris, 2003 }}
\enddedicatory

\endtopmatter

\document

\section{Introduction}

\subsection{Statement of the problem and principal results}
The question
on
the maximal number
of connected components
that
a real projective variety
of
a
given (multi-)degree
may have
remains one of the most difficult and
least understood
problems in
topology of real algebraic varieties.
Besides the trivial case of
varieties
of dimension zero,
essentially the only general situation where this problem
is solved is that of curves: the answer is
given by the famous Harnack inequality in the case of plane curves
\cite{Harnack}, and by a combination of
the
Castelnuovo--Halphen
\cite{Castel}, \cite{Halphen} and Harnack--Klein \cite{Klein}
inequalities in the case of curves in
projective spaces
of higher dimension,
see~\cite{Hilbertpaper} and~\cite{Pecker}.

The immediate
generalization of the Harnack inequality given by the
Smith theory, the so called \emph{Smith inequality}
(see, \eg,~\cite{DK}),
involves all
Betti numbers of the real
part
and the resulting bound is too rough
when applied
to the problem of the number of connected components in a
straightforward manner (see, \eg, the discussion in
Section~\ref{s.maximality}).

In
this
paper,
we address
the problem of
the maximal
number of connected components in the case of varieties defined by
equations of degree two, \ie, complete intersections of quadrics.
To be more precise, let us
denote by
$$
\MAX(r,N),\quad 0\le r\le N-1,
$$
the maximal number of connected components
that
a regular
complete intersection of $r+1$
real
quadrics in $\RpN$ can
have.
Certainly,
as we study regular complete intersections of even degree,
the actual number of connected components covers
the whole range of  values between~$0$ and $\MAX(r,N)$.

In the following
three
extremal cases,
the answer
is easy and well known:
\Dashes
\dash
$\MAX(0,N)=1$ for all $N\ge2$ \rom(a single quadric\rom),
\dash
$\MAX(1,N)=2$ for all $N\ge3$ \rom(intersection of two
quadrics\rom), and
\dash
$\MAX(N-1,N)=2^N$ for all $N\ge1$
\rom(intersection of dimension zero\rom).
\endDashes

To our knowledge, very little was known in the next case $r=2$
(intersection of three quadrics); even the fact that
$\MAX(2,N)\to\infty$ as $N\to\infty$
does not seem to have been
observed before.
Our principal result here is the following theorem, providing
a lower bound $\frac14N^2+O(N)$ and an upper bound
$\frac38N^2+O(N)$
for $\MAX(2,N)$.

\theorem\label{th.B.2}
For all $N\ge4$, one has
$$
\frac14(N-1)(N+5)-2<\MAX(2,N)\le\frac32k(k-1)+2,
$$
where $k=[\frac12N]+1$.
\endtheorem

The proof
of Theorem~\ref{th.B.2} found in Section~\ref{proof.B.2}
is based on a real version of
the Dixon
correspondence~\cite{Dixon}
between nets of quadrics
(\ie, linear systems generated by three independent quadrics) and
plane curves equipped with a non-vanishing
even
theta characteristic.
Another tool is
a spectral sequence due to Agrachev~\cite{AgraNote}, which
computes
the
homology of a complete intersection of quadrics in terms of its
spectral variety. The following intermediate result seems
to be of an independent interest.

\definition\label{def.Hilbert}
Define the \emph{Hilbert number} $\omax{d}$ as the maximal number
of ovals of the submaximal depth $[d/2]-1$ that a nonsingular real
plane algebraic curve of degree~$d$ may have. (Recall that the
depth of an oval of a curve of degree $d=N+1$ does not exceed
$[(N+1)/2]$. A brief introduction to topology of nonsingular real
plane algebraic curves is found in Section~\ref{s.topology}.)
\enddefinition

\theorem\label{th.estimate}
For any integer $N\ge4$, one has
$$
\omax{N+1}\le\MAX(2,N)\le\omax{N+1}+1.
$$
\endtheorem

This theorem is proved in Section~\ref{proof.estimate}. The few
known values of $\omax{N+1}$ and $\MAX(2,N)$ are given by the
following table.
$$
\hbox{\hss\vbox{\offinterlineskip\halign{\vrule height10.5pt depth3.5pt
 \quad\hss\strut$#$\quad\hss\vrule&&\quad\hss$#$\quad\hss\vrule\cr\noalign{\hrule}
N&3&4&5&6&7\cr\noalign{\hrule}
\omax{N+1}&4&6&9&13&17 \text{ or } 18\cr\noalign{\hrule}
\MAX(2,N)&8&6&10&
13 \text{ or } 14&17, 18, \text{ or } 19\cr\noalign{\hrule}
}}\hss}
$$

It is worth mentioning that $\MAX(2,4)=\omax5$, whereas
$\MAX(2,5)=\omax6+1$, see Sections~\ref{s.Rp4} and~\ref{s.Rp5},
respectively. (The case $N=3$ is not covered
by Theorem \ref{th.estimate}.) At present, we do not know the
precise relation between the two
sequences.

For completeness, we also discuss another extremal case, namely
that of curves. Here, the maximal number of components is attained
on the $M$-curves, and the statement should be a special case of
the general Viro--Itenberg construction
producing
maximal complete intersections of any multidegree (see \cite{IV}
for a simplified version of this construction). The result is
the following theorem, which is proved in~\S\ref{S.curves}
by means of a Harnack like construction.

\theorem\label{th.B.-2}
For all $N\ge2$, one has $\MAX(N-2,N)=2^{N-2}(N-3)+2$.
\endtheorem

\subsection{Conventions}
Unless indicated explicitly, the coefficients of
all homology and cohomology groups are~$\ZZ$.
For a compact complex curve~$C$, we freely
identify $H^1(C;R)=H_1(C;R)$ (for any coefficient ring~$R$)
\via\ the Poincar\'{e} duality.
We do not distinguish between line bundles,
invertible sheaves, and classes of linear equivalence of divisors,
switching freely from one to another.

A \emph{quod erat demonstrandum} symbol `$\square$'
after a statement means that no proof will
follow: either the statement is obvious and the proof is
straightforward, or the proof has already been explained, or a
reference is given at the beginning of the section.

\subsection{Content of the paper}
The bulk of the paper, except \S\ref{S.curves} where Theorem~\ref{th.B.-2} is proved,
is devoted to the proof of Theorems~\ref{th.B.2} and \ref{th.estimate}.
In \S\ref{theta},
we
collect the necessary
material
on the
theta characteristics, including
the real
version of the theory.
In \S\ref{systems}, we introduce and study the spectral curve of a
net and discuss Dixon's correspondence. The
aim
is to introduce
the
$\Spin$- and index (semi-)orientations
of the real part of the spectral curve,
the former
coming from the theta characteristic
and the latter, directly from the topology of the net,
and to show that
the two semi-orientations
coincide. In
\S\ref{topology.zerolocus}, we introduce
the Agrachev
spectral sequence
computing the Betti
numbers
of the common zero locus of a net in terms of
its index function.
The sequence is used to prove
Theorem~\ref{th.estimate},
relating the number $\MAX(2,N)$ in question and topology of real
plane algebraic curves of of degree $N+1$.
Then, we cite a few known estimates on the number of the ovals of
a curve (see Corollaries~\ref{Petrovsky} and~\ref{Hilbert}) and
deduce
Theorem~\ref{th.B.2}.
Finally,
in \S\ref{S.concluding} we discuss
a few
particular cases
of nets
and address
several
related questions.

\subsection{Acknowledgements}
This paper
was originally inspired
by
the following question
suggested
to us by
D.~Pasechnik and B.~Shapiro: is the number of connected
components of an intersection of $(r+1)$ real quadrics in~$\RpN$
bounded by a constant $C(r)$ independent of~$N$?
The paper was essentially completed during our stay at
\emph{Centre Interfacultaire Bernoulli}, \emph{\'{E}cole
polytechnique f\'{e}d\'{e}rale de Lausanne},
and the final
version was prepared
during the stay of the second and third authors at
the
\emph{Max-Planck-Institut f\"{u}r Mathematik}, Bonn.
We
are
grateful
to
these institutions
for their hospitality and excellent
working conditions.

\section{Theta characteristics}\label{theta}

For the reader's convenience,
we cite a number of known results related to the (real)
theta characteristics on algebraic curves. Appropriate
references are given at the beginning of each section.

To avoid various `boundary effects', {\it we only consider curves of
genus at least~$2$. For real curves, we
assume that the real part is nonempty.}

\subsection{Complex curves \rm(see \cite{Atiyah}, \cite{Mumford})}
Recall that a \emph{theta characteristic} on a nonsingular compact
complex curve~$C$
is a line bundle~$\theta$
on~$C$ such that $\theta^2$ is isomorphic to
the canonical bundle~$K_C$.
In
topological terms, a theta characteristic is merely a $\Spin$-structure
on the topological surface~$C$. One
associates with a theta characteristic~$\theta$ the
integer $h(\theta)=\dim H^0(C;\theta)$ and its
$\ZZ$-residue $\phi(\theta)=h(\theta)\bmod2$.

Let $\sset\subset\Pic^{g-1}C$ be the set of theta characteristics
on~$C$.
The map $\phi\:\sset\to \ZZ$ has the following fundamental
properties:
\roster
\item\local{phi.1}
$\phi$ is preserved under deformations;
\item\local{phi.2}
$\phi$ is a
quadratic extension
of the intersection index form;
\item\local{phi.3}
 $\Arf\phi=0$ (equivalently,
$\phi$ vanishes at $2^{g-1}(2^g+1)$ points).
\endroster
To precise the meaning of items~\loccit{phi.2} and~\loccit{phi.3},
notice that $\sset$ is an affine space over $H^1(C)$,
so that \loccit{phi.2} is equivalent to the identity
$$
\phi(a+x+y)-\phi(a+x)-\phi(a+y)+\phi(a)=\langle x,y\rangle,
$$
while $\Arf\phi$ is the usual $\Arf$-invariant of~$\phi$ after
$\sset$ is identified with $H^1(C)$ by choosing for zero any
element $\theta\in\sset$ with $\phi(\theta)=0$.

A theta characteristic $\theta$ is called \emph{even} if
$\phi(\theta)=0$; otherwise it is called \emph {odd}. An even
theta characteristic is called
\emph{non-vanishing} (or \emph{non-zero})
if
$h(\theta)=0$.

Recall that there are canonical bijections between the set of
$\Spin$-structures on~$C$, the set of quadratic extensions of
the intersection index form on $H^1(C)$, and the set of
theta characteristics on~$C$. In particular, a
theta characteristic $\theta\in\sset$ is uniquely determined by the
quadratic function $\phi_\theta$ on $H^1(C)$ given by $\phi_\theta(x)
=\phi(\theta+x)-\phi(\theta)$. One has
$\Arf\phi_\theta=\phi(\theta)$.

\paragraph\label{Rokhlin.form}
The
moduli space of pairs
$(C,\theta)$, where $C$ is a curve of a given genus and $\theta$ is a
theta characteristic on~$C$, has two connected components,
formed by even and odd theta characteristics.
If $C$
is restricted to nonsingular plane curves of a
given degree~$d$, the result is almost the same. Namely, if $d$ is
even, there are still two connected components, while if $d$ is
odd, there
is an additional component formed by the pairs
$(C,\frac12(d-3)H)$,
where $H$ is the hyperplane section divisor.
In topological terms, the extra component consists of the pairs
$(C,\Cal R)$, where $\Cal R$ is the Rokhlin function (see \cite{Rokhlin});
it is even if $d=\pm1\bmod8$ and odd if $d=\pm3\bmod8$.
The other theta characteristics
still form two connected components,
distinguished by the parity.

\subsection{Real curves \rm(see \cite{DN}, \cite{GH})}\label{orbits}
Now, let~$C$ be a real curve, \ie, a complex curve equipped with an
anti-holomorphic involution $c\:C\to C$ (a \emph{real structure}).
Recall that we always assume that the genus $g=g(C)>1$
and that
the \emph{real part} $\CR=\Fix c$ is nonempty.

Consider
the set
$$
\sset_\R=\sset\cap \Pic^{g-1}_\R C
$$
of \emph{real} (\ie, $c$-invariant)
theta characteristics.
Under the assumptions above,
there are canonical bijections between $\sset_\R$, the set of
$c$-invariant $\Spin$-structures on~$C$, and
the set of $c_*$-invariant quadratic extensions of
the intersection index form on
the group
$H^1(C)$.

The set
$\sset_\R$
of real theta characteristics
is a principal homogeneous space over
the $\ZZ$-torus $(J_\R)_2$
of torsion~$2$ elements in the real part $J_\R(C)$ of
the Jacobian~$J(C)$. In particular, $\Card\sset_\R=2^{g+r}$,
where $r=b_0(\CR)-1$.
More precisely, $\sset_\R$ admits a free action
of the $\ZZ$-torus $(J_\R^0)_2$, where $J_\R^0\subset J_\R$ is
the component of zero. This action has $2^r$ orbits,
which are distinguished by the restrictions
$\phi_\theta\:(J_\R^0)_2\to\ZZ$, which are linear forms.
Indeed, the
form~$\phi_\theta$ depends only on the orbit of
an element $\theta\in\sset_\R$, and the forms defined by
elements~$\theta_1$, $\theta_2$ in distinct orbits of the action
differ.

Alternatively, one can distinguish the orbits above as follows.
Realize an element $\theta\in\sset_\R$
by a real divisor~$D$ and, for each real component
$C_i\subset\CR, 1\le i\le r+1$, count the residue
$c_i(\theta)=\Card(C_i\cap D)\bmod2$.
The residues $(c_i(\theta))\in\ZZ^{r+1}$ are subject to relation
$\sum c_i(\theta)=g-1\mod 2$ and determine the orbit.

In most cases,
within each of the above $2^r$
orbits, the numbers of even and odd theta characteristics
coincide. The only exception to this rule is the orbit given by
$c_1(\theta)=\ldots=c_r(\theta)=1$ in the case when $C$ is a dividing curve.

\lemma\label{stability}
With one exception,
any \rom(real\rom) even theta characteristic on a
\rom(real\rom) nonsingular plane curve
becomes
non-vanishing
after a small \rom(real\rom) perturbation of the
curve in the plane.
The exception is Rokhlin's
theta characteristic $\frac12(d-3)H$ on a curve of degree
$d=\pm1\bmod8$, see~\ref{Rokhlin.form}.
\endlemma

\proof As is well known, the vanishing of a theta characteristic
is
an
analytic
condition
with respect to the coefficients of the
curve. (Essentially, this statement follows from the fact that the
Riemann $\Theta$-divisor depends on the coefficients
analytically.) Hence, in the space of pairs $(C,\theta)$, where
$C$ is a
nonsingular plane curve of degree $d$ and $\theta$ is
an even
theta
characteristic on~$C$, the pairs $(C,\theta)$
with non-vanishing $\theta$ 
form a Zariski open set. Since there does exist a
curve of degree $d$ with a
non-vanishing even
theta characteristic (\eg, any nonsingular
spectral curve, see \ref{nonvanishing}), this set is nonempty and
hence dense in the (only) component formed by the even theta
characteristics other than $\frac12(d-3)H$.
\endproof

\section{Linear systems of quadrics}\label{systems}

\subsection{Preliminaries}\label{s.notation}
Consider an injective linear map $x\mapsto\tq_x$ from $\C^{r+1}$
to the space $S^2\C^{N+1}$ of homogeneous quadratic polynomials
on~$\C^{N+1}$. It defines a linear system of quadrics in $\CpN$ of
dimension~$r$, \ie, an $r$-subspace in the projective  space
$\CC_2(\CpN)$ of quadrics. Conversely, any linear system of
quadrics is defined by a unique, up to obvious equivalence, linear
map as above.

Occasionally, we will fix coordinates $(u_0,\ldots,u_N)$ in
$\C^{N+1}$ and represent $q_x$ by a matrix $Q_x$, so that
$q_x(u)=\langle Q_xu,u\rangle$. Clearly, the map $x\mapsto Q_x$ is
also linear.

Define the common zero set
$$
V=\bigl\{u\in\CpN\bigm|
 \text{$q_x(u)=0$ for all $x\in\C^{r+1}$}\bigr\}\subset\CpN
$$
and the \emph{Lagrange hypersurface}
$$
L=\bigl\{
(x,u)\in\Cpr\times\CpN\bigm|
 q_x(u)=0\bigr\}\subset\Cpr\times\CpN.
$$
(As usual, the vanishing condition $q_x(u)=0$ does not depend on
the choice of the representatives of~$x$ and~$u$.) The following
statement is straightforward.

\lemma\label{Lregular}
An intersection of quadrics~$V$ is regular if and only if the
associated Lagrange hypersurface~$L$ is nonsingular.
\qed
\endlemma

\subsection{The spectral variety}\label{s.spectral}
Define the \emph{spectral variety}~$C$ of a linear system of
quadrics $x\mapsto\tq_x$ \via\
$$
C=\bigl\{x\in\Cpr\bigm|\det Q_x=0\bigr\}\subset\Cpr.
$$
Clearly, this definition does not depend on the choice of
the matrix representation $x\mapsto Q_x$: the spectral variety is
formed by the elements of the linear system which are singular
quadrics. More precisely (as a scheme),
$C$ is the intersection of the linear
system with the discriminant hypersurface
$\Delta\subset\CC_2(\CpN)$.

In what follows, we
assume that $C$ is a proper subset of $\Cpr$, \ie, we exclude the possibility $C=\Cpr$,
as in this
case all quadrics in the system have a common singular point and,
hence, $V$ is not a regular intersection. Under this assumption,
$C\subset\Cpr$ is
a hypersurface of degree~$N+1$, possibly not reduced. By the
dimension argument, $C$ is necessarily singular whenever $r\ge3$.
Furthermore, even in the case $r=1$ or~$2$ ({\it pencils} or {\it nets}), the
spectral variety of a regular intersection may still be
singular.

\lemma\label{pencil}
Let $x$ be an isolated point of the spectral variety~$C$ of a
pencil. Then $x$ is a simple point of~$C$ if and only if the
quadric $\{q_x=0\}$ has
a single singular point, and this point is not
a base point of the pencil.
\qed
\endlemma

\corollary\label{Cregularity}
For any $r$,
if $x\in C$ is a smooth point, then $\corank\tq_x=1$,
\ie, the quadric $\{\tq_x(u)=0\}$ has
a single
singular point.
\endcorollary

\proof
Restrict the system to a generic pencil through the point.
\endproof

\lemma\label{CtoVregularity}
If the spectral
variety~$C$ of a linear system is nonsingular,
then the complete intersection $V$ is regular.
\endlemma

\proof
It is easy to see that the common zero set~$V$ of a linear system
is a regular complete intersection if and only if none of the
members of the system has a singular point in~$V$. Thus, it
suffices to observe that a generic pencil through a point $x\in C$
is transversal to~$C$; hence, $x$ is a simple point of its
discriminant variety and the only singular point of $\{q_x=0\}$ is
not in~$V$, see Lemma~\ref{pencil}.
\endproof

\subsection{The Dixon construction
\rm(see \cite{Dixon}, \cite{Dolg})}\label{s.Dixon} From now on, we
confine ourselves to the case of nets, \ie, $r=2$. As explained
in~\ref{s.spectral}, each net gives rise to its spectral curve,
which is a curve
$C \subset \Cp2$
of degree $d=N+1$;
if the net is
generic, $C$ is nonsingular.

\paragraph\label{nonvanishing}
Assume that the spectral curve
$C$ is nonsingular. Then, at each point $x\in C$, the
kernel $\Ker Q_x\subset\R^{N+1}$ is a $1$-subspace, see
Lemma~\ref{Cregularity}. The correspondence $x\mapsto\Ker Q_x$
defines a line bundle $\CK$ on~$C$ or, after a twist, a line
bundle $\CL=\CK(d-1)$.
The latter has the following properties: $\CL^2=\CO_C(d-1)$
(so that $\deg\CL=\frac12 d(d-1)$\,) and $H^0(C,\CL(-1))=0$. Thus,
switching to $\theta=\CL(-1)$, we obtain a non-vanishing
even
theta
characteristic on~$C$; it is called the
\emph{spectral theta characteristic} of the net.

The following theorem is due to A.~C.~Dixon~\cite{Dixon}.

\theorem\label{DixonTh}
Given a
non-vanishing
even
theta characteristic~$\theta$
on a nonsingular plane curve~$C$ of degree $N+1$,
there exists
a unique, up to projective
transformation of $\CpN$,
net of quadrics in $\CpN$
such that
$C$ is its spectral curve and $\theta$ is its spectral
theta characteristic.
\qed
\endtheorem

\paragraph\label{s.construction}
The original proof by Dixon contains an explicit construction of
the net.
We outline this construction below.
Pick a basis
$\phi_{11},\phi_{12},\dots, \phi_{1d}\in H^0(C,\CL)$ and
let
$$
v_{11}=\phi_{11}^2,\
v_{12}=\phi_{11}\phi_{12},\ \dots,\
v_{1d}=\phi_{11}\phi_{1d}\in H^0(C, \CL^2).
$$
Since the restriction map
$H^0(\Cp2;\CO_{\Cp2}(d-1))\to H^0(C;\CO_C(d-1))=H^0(C;\CL^2)$ is
onto, we can regard~$v_{1i}$
as
homogeneous polynomials of degree $d-1$
in
the coordinates $x_0,x_1,x_2$ in~$\Cp2$.
Let also $U(x_0,x_1,x_2)=0$ be the equation of~$C$.
The curve
$\{v_{12}=0\}$
passes through
all points of intersection of~$C$ and $\{v_{11} = 0\}$.
Hence, there
are
homogeneous polynomials $v_{22}$, $w_{1122}$ of degrees $d-1$,
$d-2$, respectively, such that
$$
v_{12}^2=v_{11}v_{22}-Uw_{1122}.
$$
In the same way, we get polynomials $v_{rs}$, $w_{11rs}$,
$2\le s,r\le d$, such that
$$
v_{1r}v_{1s}=v_{11}v_{rs}-Uw_{11rs}.
$$
Obviously, $v_{rs}=v_{sr}$ and $w_{11rs}=w_{11sr}$.
It is shown in~\cite{Dixon} that
\roster
\item
the algebraic complement $A_{rs}$ in the
$(d\times d)$
symmetric matrix $[v_{ij}]$
is of the form $U^{d-2}\beta_{rs}$,
where $\beta_{rs}$ are certain linear forms,
and
\item
the determinant $\det[\beta_{ij}]$
is a constant non-zero multiple of~$U$.
\endroster
(It is the
non-vanishing of the theta characteristic that is used to show
that the latter determinant is non-zero.) Thus, $C$ is the
spectral curve of the net $Q_x=[\beta_{ij}]$.


It is immediate that the construction works over any field of
characteristic zero. Hence, we obtain the following real version
of the Dixon theorem.

\theorem\label{Dixon.real}
Given
a nonsingular real plane curve~$C$ of
degree $N+1\ge 4$
with nonempty real part,
and a
real
non-vanishing
even
theta characteristic~$\theta$
on $C$, there exists a unique, up to real projective
transformation of $\RpN$, real net of quadrics in $\RpN$ such that
$C$ is its spectral curve and $\theta$ is its spectral theta
characteristic.
\qed
\endtheorem

\subsection{The $\Spin$-orientation \rm({\it cf}\. \cite{Natanzon},
\cite{NataBook})}\label{loops}
Let $(C,c)$ be a real curve equipped with a real
theta characteristic~$\theta$. As above, assume
that $\CR\ne\emptyset$.
Then, the
real structure of~$C$ lifts to a real
structure (\ie, a fiberwise anti-linear involution)
$c\:\theta\to\theta$, which is unique up to
a phase
factor
$e^{i\phi}$, $\phi\in\R$.
If an isomorphism $\theta^2=K_C$ is fixed, one can choose a lift
compatible with the canonical action of~$c$ on~$K_C$; such a lift is
unique up to multiplication by~$i$.

Fix a lift $c\:\theta\to\theta$ as above and pick a $c$-real meromorphic
section~$\omega$ of~$\theta$.
Then, $\omega^2$ is a real
meromorphic $1$-form with zeros and poles
of even multiplicities. Therefore, it determines an orientation of
$\CR$.
This orientation does not depend on the choice of~$\omega$,
and it is reversed when switching from~$c$ to~$ic$.
Thus, it is, in fact, a semi-orientation of~$\CR$; it is called
the \emph{$\Spin$-orientation} defined by~$\theta$.

The definition above can be made closer to original
Dixon's construction outlined in~\ref{s.construction}. One can
replace~$\omega$ by a meromorphic
section~$\omega'$
of $\theta(1)$
and treat $(\omega')^2$ as a real meromorphic $1$-form
with values in $\CO_C(2)$; the latter is trivial over~$\CR$.

The
following, more topological, definition is equivalent to
the previous one.
Recall that a semi-orientation is essentially a rule comparing
orientations of pairs of components.
Let~$\theta$ be a $c$-invariant
$\Spin$-structure on~$C$, and let $\phi_\theta\:H_1(C)\to\ZZ$ be
the associated quadratic extension.
Pick a point~$p_i$
on each real component~$C_i$ of~$\CR$. For each pair
$p_i$, $p_j$, $i\ne j$, pick a simple smooth
path connecting~$p_i$ and~$p_j$
in the complement $C\sminus\CR$ and transversal to~$\CR$ at the ends,
and let~$\gamma_{ij}$ be
the loop obtained by
combining the path with its $c$-conjugate. Pick an orientation
of~$C_i$ and transfer it to~$C_j$ by a vector field normal
to the path above. The two orientations are considered
coherent with respect
to the $\Spin$-orientation defined by~$\theta$
if and only if
$\phi_\theta([\gamma_{ij}])=0$.



\lemma\label{arbitrarySpin}
Assuming that $\CR\ne\emptyset$,
any semi-orientation of~$\CR$
is the $\Spin$-orientation for a suitable
real even theta characteristic.
\qed
\endlemma

\proof
Observe that
the $c_*$-invariant classes
$\{[\gamma_{1i}],[C_i]\}$ with $i\ge2$ (see above)
form a standard symplectic basis in a
certain nondegenerate
$c_*$-invariant subgroup $S\subset H_1(C)$.
On the complement~$S^\perp$,
one can pick any $c_*$-invariant quadratic
extension with
$\Arf$-invariant~$0$.
Then,
the values $\phi_\theta([\gamma_{1i}])$ can be chosen arbitrarily
(thus producing any given $\Spin$-orientation), and
the values
$\phi_\theta([C_i])$, $i\ge2$, can be adjusted
(\eg, made all~$0$)
to make the resulting theta characteristic even.
\endproof

\paragraph\label{alt.orientation}
If $d=\pm1\bmod8$ and $\theta$ is the exceptional theta characteristic
$\frac12(d-3)H$, see \ref{Rokhlin.form},
the $\Spin$-orientation defined by~$\theta$ is
given by the residue
$\res(p^2\Omega/U)$, where $U=0$ is the equation of~$C$ as above,
$\Omega=x_0dx_1\wedge dx_2-x_1dx_0\wedge dx_2+x_2dx_0\wedge dx_1$
is a non-vanishing section of
$K_{\Cp2}(3)\cong\CO_{\Cp2}$,
and $p$ is any
real homogeneous polynomial of degree $(d-3)/2$. In affine
coordinates, this orientation is given by $p^2\,dx\wedge dy/dU$.
Such a
semi-orientation is called \emph{alternating}: it is the only
semi-orientation of $\CR$ induced by alternating orientations of
the components of $\Rp{2}\sminus \CR$.

\subsection{The index function
\rm({\it cf}\.~\cite{Agrachev})}\label{s.index}
Fix a real dimension~$r$
linear system $x\mapsto\tq_x$ of quadrics in $\CpN$. Consider the
sphere
$S^r=(\R^{r+1}\sminus0)\!/\R_+$
and denote by $\tC\subset S^r$ the pull-back of
the real part $\CR\subset\Rpr$ of
the spectral variety under the double covering $S^r\to\Rpr$.
Define the \emph{index function}
$$
\ind\:S^{r}\to\Z
$$
by sending a point $x\in S^r$ to
the negative index of inertia of the quadratic form~$\tq_x$.
The following statement is obvious.
(For
item~\iref{prop.ind}{ind.jump}, one should use
Corollary~\ref{Cregularity}.)

\proposition\label{prop.ind}
The index function~$\ind$ has the following properties\rom:
\roster
\item\local{ind.scont}
$\ind$ is lower semicontinuous\rom;
\item\local{ind.const}
$\ind$ is locally constant on $S^r\sminus\tC$\rom;
\item\local{ind.jump}
$\ind$ jumps by~$\pm1$
when crossing~$\tC$ transversally at its regular point\rom;
\item\local{ind.-x}
one has $\ind(-x)=N+1-(\ind x+\corank\tq_x)$.
\qed
\endroster
\endproposition

Due to~\iref{prop.ind}{ind.jump}, $\ind$
defines a coorientation of~$\tC$ at all its smooth points. This
coorientation is reversed by the antipodal map $a\:S^r\to S^r$,
$x\mapsto-x$, see~\iref{prop.ind}{ind.-x}.
If $r=2$ and the real part~$\CR$ of the spectral curve is nonsingular, this
coorientation defines a semi-orientation of~$\CR$ as follows: pick
an orientation of~$S^2$ and use it to convert the coorientation to
an orientation~$\tilde o$ of~$\tC$;
since the
antipodal map~$a$ is orientation
reversing, $\tilde o$ is preserved by~$a$ and, hence, descends to an
orientation~$o$ of~$\CR$. The latter is defined up to the total
reversing (due to the initial choice of an orientation of~$S^2$);
hence, it is in fact a semi-orientation. It is called the
\emph{index orientation} of~$\CR$.

Conversely, any semi-orientation of~$\CR$ can be defined as above
by a function $\ind\:S^r\to\Z$
satisfying~\iref{prop.ind}{ind.scont}--\ditto{ind.-x}; the latter
is unique up to the antipodal map.


\theorem[Theorem \rm(\cf.~\cite{Vinnikov})]\label{index-to-spin}
Assume that the real part~$\CR$ of the spectral curve of a
real net of quadrics is
nonsingular. Then the
index orientation of~$\CR$ coincides with its
$\Spin$-orientation defined by the spectral theta characteristic.
\endtheorem

\proof
The
$\Spin$-semiorientation of $\CR$ is given by the residue
$\res(v_{11}\Omega/U)$, \cf.~\ref{alt.orientation}, and
it is sufficient to check that the index function is
larger on the side of $U=0$ where $v_{11}/U>0$. Since
the kernel $\sum x_iQ_i$, regarded as a section of the
projectivization of the trivial bundle over~$C$,
is given by $v=(v_{11},v_{12},\dots, v_{1d})$, it remains to
notice that, over~$\CR$, one has
$\sum x_i\langle Q_iv, v\rangle =v_{11}\det(\sum x_iQ_i)= v_{11} U$.
\endproof

\theorem\label{main}
Let $C$ be a nonsingular real plane curve of degree $d=N+1$
and
with nonempty real part,
and let~$o$ be a
semi-orientation of~$\CR$. Assume that either $d\ne\pm1\bmod8$ or
$o$ is not the alternating semi-orientation,
see~\ref{alt.orientation}.
Then,
after a small real
perturbation of~$C$, there exists
a regular intersection of three real quadrics in~$\RpN$
which has $C$ as its spectral curve and $o$ as its
spectral $\Spin$-orientation.
\endtheorem

\Remark\label{rem.index}
According to~\ref{alt.orientation}, the only case not covered by
Theorem~\ref{main} is when $d=\pm1\bmod8$ and the index
function~$\ind$
assumes only the two middle
values $(d\pm1)/2$.
\endRemark

\proof[Proof of~\ref{main}] By Lemma~\ref{arbitrarySpin}, there
exists a real even theta characteristic $\theta$ that has~$o$ as
its $\Spin$-orientation. Using Lemma \ref{stability}, one can
make~$\theta$ non-vanishing by a small real perturbation of~$C$,
and it remains to apply Theorem~\ref{Dixon.real}.
\endproof

\section{The topology of the zero locus of a net}\label{topology.zerolocus}

\subsection{The spectral sequence}
Consider a real dimension~$r$ linear system $x\mapsto\tq_x$ of
quadrics in~$\CpN$, see Section~\ref{s.notation} for the notation.
Let $\VR\subset\RpN$, $\CR\subset\Rpr$, and
$\LR\subset\Rpr\times\RpN$ be the real parts of the common zero
set, spectral variety, and Lagrange hypersurface, respectively.

Consider the sphere $S^r=(\R^{r+1}\sminus0)\!/\R_+$
and the lift
$\tC\subset S^r$ of~$\CR$, \cf.
\ref{s.index}, and let
$\tL=\{(x,u)\in\Sr\times\RpN\,|\,\tq_x(u)=0\}\subset\Sr\times\RpN$
be the lift of~$\LR$ and
$$
\Lm=\{(x,u)\in\Sr\times\RpN\,|\,\tq_x(u)>0\}\subset\Sr\times\RpN
$$
its \emph{positive complement}. (Clearly, the conditions
$\tq_x(u)=0$ and $\tq_x(u)>0$ do not depend on the choice
of
representatives of $x\in S^r$ and $u\in\RpN$.)

\lemma\label{L+=complement}
The projection $\Sr\times\RpN\to\RpN$
restricts to a homotopy equivalence $\Lm\to\RpN\sminus\VR$.
\endlemma

\proof
Denote by~$p$ the restriction of the projection to~$\Lm$.
The pull-back $p^{-1}(u)$ of a point $u\in\VR$ is empty; hence,
$p$ sends~$\Lm$ to $\RpN\sminus\VR$. On the other hand, the
restriction $\Lm\to\RpN\sminus\VR$ is a locally trivial fibration,
and for each
point $u\in\RpN\sminus\VR$, the fiber $p^{-1}(u)$ is the
open hemisphere
$\{x\in\Sr\,|\,q_x(u)>0\}$, hence
contractible.
\endproof

%

\proposition\label{b(L)}
Let $r=2$. Then, one has $b^0(\Lm)=b^1(\Lm)=1$,
and if $\VR$ is nonsingular,
also $b^2(\Lm)=b^0(\VR)+1$.
\endproposition

\proof Due to Lemma~\ref{L+=complement}, one has
$H^*(\Lm)=H^*(\RpN\sminus\VR)$, and the statement follows from the
Poincar\'{e}-Lefschetz duality
$H^i(\RpN\sminus\VR)=H_{N-i}(\RpN,\VR)$ and the exact sequence of
the pair $(\RpN,\VR)$. For the last statement one needs, in
addition, to know that the inclusion homomorphism
$H_{N-3}(\VR)\to H_{N-3}(\RpN)$ is trivial, \ie, that every
$3$-plane~$P$ intersects each component of~$\VR$ at an even number of
points. By restricting the system to a $4$-plane containing~$P$,
one reduces the problem to the case $N=4$. In this case,
$V\subset\CpN$ is the canonical embedding of a genus~$5$
curve, \cf. Section~\ref{s.Rp4},
and the statement is obvious.
\endproof

 From now on, we assume that the real part~$\CR$ of the spectral
hypersurface is nonsingular.
Consider the ascending filtration
$$
\varnothing=\SS_{-1}\subset \SS_0\subset \SS_1\subset\ldots
 \subset \SS_{N+1}=\Sr,\quad
\SS_i=\{x\in\Sr\,|\,\ind x\le i\}.
\eqtag\label{filtration}
$$
Due to Proposition~\iref{prop.ind}{ind.scont},
all $\SS_i$ are closed subsets.

\theorem[Theorem \rm(\cf.~\cite{AgraNote})]\label{th.ss}
There is a spectral sequence
$$
E_2^{pq}=H^p(\SS_{N-q})\Rightarrow H^{p+q}(\Lm).
$$
\endtheorem

\proof The sequence in question is the Leray spectral sequence of
the projection $\pi\:\Lm\to\Sr$. Let~$\CZ_2$ be the constant sheaf
on~$\Lm$ with the fiber~$\ZZ$. Then the sequence is
$E_2^{pq}=H^p(\Sr,R^q\pi_*\CZ_2)\Rightarrow H^{p+q}(\Lm)$. Given a
point $x\in\Sr$, the stalk
$(R^q\pi_*\CZ_2)|_x$ equals $H^q(\pi^{-1}U_x)$, where $U_x\ni x$ is
a small neighborhood of~$x$ regular with respect to a
triangulation of~$\Sr$ compatible with the filtration. If
$x\notin\tC$ and $\ind x=i$, then
$$
\pi^{-1}U_x\sim\pi^{-1}x=\{u\in\RpN\,|\,\tq_x(u)>0\}
 \sim\Rp{N-i}.
$$
(The fiber $\pi^{-1}x$ is a $D^i$-bundle over $\Rp{N-i}$;
if $i=N+1$, the fiber is empty.)
Besides, if $x\in\tC$ and $\ind x=i$, then
$\pi^{-1}U_x\sim\pi^{-1}x'$, where
$x'\in(U_x\cap\SS_i)\sminus\tC$.
Thus,
$(R^q\pi_*\CZ_2)|_x=\ZZ$ or~$0$ if $x$ does (respectively, does
not)
belong to $\SS_{N-q}$, and the statement is immediate.
\endproof

\Remark
It is worth
mentioning that there are spectral sequences, similar to the
one introduced in
Theorem~\ref{th.ss}, which compute the cohomology of
the double coverings of $\RpN\sminus\VR$ and $\VR$ sitting in
$S^N$, see~\cite{Agrachev}.
\endRemark

\subsection{Elements of topology of real plane curves}\label{s.topology}
Let $C\subset\Cp2$ be a nonsingular real
curve of
degree~$d$. Recall that the
real part~$\CR$
splits into
a number of \emph{ovals} (\ie, embedded
circles
contractible in $\Rp2$) and, if $d$ is odd, one \emph{one-sided
component} (\ie, an embedded circle isotopic to~$\Rp1$.)
The complement
of each oval~$\oval$ has two connected
components, exactly one of them being
contractible; this contractible component is called the
\emph{interior} of~$\oval$.

On the set of ovals of~$\CR$, there is a natural partial order: an
oval~$\oval$ is said to \emph{contain}
another oval~$\oval'$, $\oval\prec\oval'$, if $\oval'$ lies in
the interior of~$\oval$. An oval is called \emph{empty}
if it does not
contain another oval.
The \emph{depth} $\dpth\oval$ of an
oval~$\oval$ is the number of elements in the maximal descending
chain starting at~$\oval$. (Such a chain is unique.)
Every oval~$\oval$ of depth $>1$ has a unique immediate predecessor;
it is denoted by $\pred\oval$.

A \emph{nest} of~$C$
is a linearly ordered chain of ovals of~$\CR$;
the \emph{depth} of a nest
is the number of its elements.
The following statement is a simple and well known
consequence of the B\'{e}zout theorem.

\proposition\label{Bezout} Let $C$ be a nonsingular real plane
curve of degree~$d$. Then \roster \item $C$ cannot have a nest of
depth
greater than
$\Dmax=\Dmax(d)=[d/2]$\rom;
\item
if $C$ has a
nest
of depth $\Dmax$ \rom(a \emph{maximal
nest}\rom), it has no other ovals\rom;
\item
if $C$ has a nest $\oval_1\prec\ldots\prec\oval_k$ of depth
$k=\Dmax-1$ \rom(a \emph{submaximal nest}\rom)
but no maximal nest,
then all ovals other than
$\oval_1,\ldots,\oval_{k-1}$ are empty.
\qed
\endroster
\endproposition

Let $\pr\:S^2\to\Rp2$ be the orientation double covering, and let
$\tC=\pr^{-1}\CR$. The pull-back of an oval~$\oval$ of~$C$
consists of two disjoint circles $\oval'$, $\oval''$; such circles
are called \emph{ovals} of~$\tC$.
The antipodal map $x\mapsto-x$ of~$S^2$ induces an involution
on the set
of ovals of~$\tC$; we denote it by bar,
$\oval\mapsto\bar\oval$.
The pull-back of the one-sided
component of~$\CR$ is connected; it is called the \emph{equator}.
The \emph{tropical components} are the components of
$S^2\sminus\tC$ whose image is the (only)
component of $\Rp2\sminus\CR$ outer to all ovals.
The \emph{interior} $\int\oval$ of an oval~$\oval$ of~$\tC$ is the
component of the complement of~$\oval$ that projects to the
interior of $\pr\oval$ in~$\Rp2$. As in the case of~$\CR$, one can
use the notion of interior to define the partial order, depth,
nests, \etc. The projection~$\pr$ and the antipodal involution
induce strictly
increasing maps of the sets of ovals.

Now, consider an (abstract) index function $\ind\:S^2\to\Z$
satisfying~\iref{prop.ind}{ind.scont}--\ditto{ind.-x} (where
$N+1=d$ and $\corank\tq_x=\chi_{\tilde C}(x)$ is the characteristic
function of~$\tC$) and use it to define the filtration $\SS_*$ as
in~\eqref{filtration}.
For an oval~$\oval$ of~$\tC$, define
$\oind\oval$ as the value of~$\ind$ immediately inside~$\oval$.
(Note that, in general, $\oind\oval$ is \emph{not} the restriction
of~$\ind$ to~$\oval$ as a subset of~$S^r$.)
Then, in view of Proposition~\ref{prop.ind}, one has
$$
\gather
\tfrac12N\le\ind|_T\le\tfrac12N+1\quad
 \text{for each tropical component~$T$},
\eqtag\label{T.ind}\\
\oind\oval\le N+1-(\Dmax-\dpth\oval)\quad
 \text{for each oval~$\oval$},
\eqtag\label{max.ind}\\
\oind\oval=N+1-(\Dmax-\dpth\oval)\bmod2\quad
 \text{if $N$ is odd}.
\eqtag\label{cong.ind}\\
\endgather
$$
Here, as above, $\Dmax=[(N+1)/2]$.
Keeping in mind the
applications, we state the restrictions in terms of $N=d-1$.
(Certainly, congruence~\eqref{cong.ind} simplifies to
$i(\oval)+\dpth\oval=\Dmax\bmod2$; however, we leave it in
the form convenient for the further applications.)

\corollary\label{tropics}
If $N\ge5$, then $\SS_{N-2}$ contains the tropical components.
\qed
\endcorollary

\lemma\label{deep.nest}
Assume that $b^0(\SS_q)>1$ for some integer $q>\frac12N$.
Then $\tC$ has a nest $\oval\prec\oval'$ such that
$\oind\oval=q+1$
and $\oind{\oval'}=q$.
\endlemma

\proof
Due to~\eqref{T.ind}, the assumption $q>\frac12N$ implies that
$\SS_q$ contains the tropical components. Then, one can take
for~$\oval'$ the oval bounding from outside another component
of~$\SS_q$, and let $\oval=\pred\oval'$.
\endproof

\lemma\label{max.nest}
Let $N\ge7$, and assume that the curve~$\tC$ has
a nest $\oval_{k-1}\prec\oval_k$, $k=\Dmax-1$, with
$\oind{\oval_{k-1}}=N-1$. Then
$\SS_{N-2}\supset S^2\sminus\int\pred\oval_{k-1}$.
\endlemma

\proof
Due to~\eqref{max.ind}, the nest $\oval_{k-1}\prec\oval_k$
in the statement
can be completed to a submaximal
nest $\oval_1\prec\ldots\prec\oval_k$.

First, assume that either $\tC$ has no maximal nest or the
innermost oval of~$\tC$ is inside~$\oval_k$.
Then
$\dpth\oval_s=s$ and $\oval_s=\pred\oval_{s+1}$ for all~$s$.
Due to~\eqref{max.ind} and~\iref{prop.ind}{ind.jump}, one has
$\oind{\oval_{k-s}}=N-s$ for $s=1,\ldots,k-1$. Then, due
to~\iref{prop.ind}{ind.-x}, $\oind{\bar\oval_{k-s}}=s+1$ for $s=1,\ldots,k-1$
and hence $\oind{\bar\oval_k}\le3$.

Proposition~\ref{Bezout} implies that all ovals other than
$\oval_{k-s}$, $\bar\oval_{k-s}$, $s=0,\ldots,k-1$, are empty. For
such an oval~$\oval$, one has $\oind\oval\le[\frac12N]+2$ if
$\dpth\oval=1$, see~\eqref{T.ind}, and $\oind\oval\le4$, $s+2$, or
$N+1-s$ if $\pred\oval=\bar\oval_k$,
$\bar\oval_{k-s}$, or $\oval_{k-s}$, respectively,
$s=1,\ldots,k-1$. From the assumption $N\ge7$, it follows that,
for any oval~$\oval$, one has
$\oind\oval\le N-2$ unless
$\oval\succcurlyeq\oval_{k-2}$. Together with
Corollary~\ref{tropics}, this observation implies the statement.

The case when~$\tC$ has another oval $\oval'\prec\oval_{i}$
for some $i\le k-1$ is treated similarly. In this case, $N=2k$ is
even, see~\eqref{cong.ind}, and, renumbering the ovals
consecutively from~$k$ down to~$0$, one has
$\oind{\oval_{k-s}}=N-s$, $\oind{\bar\oval_{k-s}}=s+1$ for
$s=1,\ldots,k$.
\endproof

Note that, in particular, the condition $N\ge7$ in Lemma~\ref{max.nest} is
necessary for~$\oval_{k-1}$ to have a predecessor, \ie,
for the statement to make sense. The remaining two
interesting cases $N=5$ and~$6$
are treated in the next lemma.

\lemma\label{max.nest.5.6}
Let $N=5$ or~$6$, and assume that the curve~$\tC$ has a nest
$\oval_1\prec\oval_2$ with $\oind{\oval_1}=N-1$. Then either
\roster
\item
for each pair $\oval$, $\bar\oval$ of
nonempty
antipodal
ovals of depth~$1$,
the set $\SS_{N-2}$ contains the interior of exactly
one of them, or
\item
$N=5$ and $\tC$ has another nested oval $\oval_3\succ\oval_2$ with
$\oind{\oval_3}=2$. \rom(These data determine the index function
uniquely.\rom)
\endroster
\endlemma

\proof
The proof repeats literally that of Lemma~\ref{max.nest}, with a
careful analysis of the inequalities that do not hold for small
values of~$N$.
\endproof

\subsection{The estimates}\label{s.estimates}
In this section, we consider a net of quadrics ($r=2$) and assume
that the spectral curve $\tC\subset S^2$ is nonsingular.
Furthermore, we can assume that the index function takes values
between $1$ and~$N$, as otherwise the net would contain an empty
quadric and one would have $V=\varnothing$. Thus, one has
$\varnothing=\SS_{-1}=\SS_0$ and $\SS_N=\SS_{N+1}=S^2$.

Denote $\imax=\max_{x\in S^2}\ind x$. Thus, we assume that
$\imax\le N$.

In addition, we can assume that $N\ge3$, as for $N\le2$ a regular
intersection of three quadrics in $\CpN$ is empty.

The spectral sequence~$E_r^{pq}$ given by Theorem~\ref{th.ss}
is concentrated in the strip $0\le p\le2$, and all potentially
nontrivial differentials are
$d_2^{0,q}\:E_2^{0,q}\to E_2^{2,q-1}$, $q\ge1$.
Furthermore, one has
$$
\gather
E_2^{0,q}=E_2^{2,q}=\ZZ,\quad E_2^{1,q}=0
 \quad\text{for $q=0,\ldots,N-\imax$},
\eqtag\label{imax}\\\allowdisplaybreak
E_2^{0,q}=E_2^{1,q}=E_2^{2,q}=0
 \quad\text{for $q\ge\imax$},\quad\text{and}
\eqtag\label{all0}\\\allowdisplaybreak
E_2^{2,q}=0\quad\text{for $q>N-\imax$}.
\eqtag\label{qmax}
\endgather
$$
In particular, it follows that $d_2^{0,q}=0$ for $q>N+1-\imax$.

\corollary\label{imax<N-1}
If $\imax\le N-2$, then $b^0(\VR)\le1$.
\qed
\endcorollary

The assertion of Corollary~\ref{imax<N-1} was first observed by
Agrachev~\cite{AgraNote}.

\lemma\label{d2}
With one exception, $d_2^{0,1}=0$.
The exception is a curve~$\tC$ with a maximal nest
$\oval_1\prec\ldots\prec\oval_k$, $k=\Dmax$, so that
$\oind{\oval_k}=N-1$ and $\oind{\oval_{k-1}}=N$.
In this exceptional case, one has $b^0(\VR)=1$.
\endlemma

\proof
Since $b^1(\Lm)=1$, see Proposition~\ref{b(L)}, and $E_2^{1,0}=0$,
the differential $d_2^{0,1}$ is nontrivial if and only if
$b^0(\SS_{N-1})=\dim E_2^{0,1}>1$. Since $N\ge3$,
the exceptional case is
covered by Lemma~\ref{deep.nest} and Proposition~\ref{Bezout}.
\endproof

\corollary\label{imax=N-1.cor}
If $\imax=N-1\ge2$, then
one has
$b^0(\SS_{N-2})-1\le b^0(\VR)\le b^0(\SS_{N-2})$.
\qed
\endcorollary

\lemma\label{imax=N-1}
Assume that $\imax=N-1\ge4$ and that
$b^0(\SS_{N-2})>1$.
Then $\beta\le b^0(\VR)\le\beta+1$,
where $\beta$ is the number of ovals of~$\CR$ of depth $\Dmax-1$.
\endlemma

\proof
In view of Corollary~\ref{imax=N-1.cor}, it suffices to show that
$b^0(\SS_{N-2})=\beta+1$. Due to Lemma~\ref{deep.nest}, the
curve~$\tC$ has a nest $\oval\prec\oval'$ with
$\oind\oval=N-1$, and then the set~$\SS_{N-2}$ is
described by Lemmas~\ref{max.nest},~\ref{max.nest.5.6} and the
assumption $\imax=N-1$. In view of Proposition~\ref{Bezout}, each
oval of~$\CR$ of depth $\dpth\oval+1$ is inside $\pr\oval$, thus
contributing an extra unit to $b^0(\SS_{N-2})$.
\endproof

\lemma\label{imax=N}
Assume that $\imax=N\ge5$. Then $b^0(\VR)$ is equal to the
number $\beta$ of ovals of~$\CR$ of depth $\Dmax-1$.
\endlemma

\proof
If $(\tC,\ind)$ is the exceptional index function mentioned in
Lemma~\ref{d2}, then $b^0(\VR)=\beta=1$, and the statement holds.
Otherwise, both differentials $d_2^{0,1}$ and $d_2^{0,2}$ vanish,
see~\eqref{qmax}, and,
using Proposition~\ref{b(L)} and~\eqref{imax}, one concludes
that
$b^0(\VR)=\dim E_2^{0,2}+\dim E_2^{1,1}=
 b^0(\SS_{N-2})+b^1(\SS_{N-1})$.

Pick an oval~$\oval'$ with $\oind{\oval'}=N$ and let
$\oval=\pred\oval'$. The topology of $\SS_{N-2}$ is given by
Lemmas~\ref{max.nest},~\ref{max.nest.5.6}. If
$\dpth\oval=\Dmax-1$, then $b^0(\VR)=\beta=1$. Otherwise
($\dpth\oval=\Dmax-2$), one has $b^0(\SS_{N-2})=\beta_-+1$ and
$b^1(\SS_{N-1})=\beta_+-1$, where $\beta_-\ge0$ and $\beta_+>0$
are the numbers of ovals $\oval''\succ\oval$ with
$\oind{\oval''}=N-2$ and~$N$, respectively; due to
Proposition~\ref{Bezout}, one has $\beta_-+\beta_+=\beta$.
\endproof

\subsection{Proof of Theorem~\ref{th.estimate}}\label{proof.estimate}
The case $N=4$ is covered by Theorem~\ref{th.B.-2}.
(Note that
$\MAX(2,4)=\omax{5}$, see Section~\ref{s.Rp4}.)
Alternatively, one can treat
this case manually, trying various index functions on a curve of
degree~$5$.

Assume that $N\ge5$. The upper bound on~$\MAX(2,N)$ follows from
Corollary~\ref{imax<N-1} and Lemmas~\ref{imax=N-1}
and~\ref{imax=N}. For the lower bound, pick a generic real
curve~$C$ of degree $d=N+1$ with $\omax{d}$ ovals of depth
$[d/2]-1$.
Select an oval~$\oval_i$ in each pair $(\oval_i,\bar\oval_i)$ of
antipodal outermost ovals of~$\tC$; if $d$ is even, make sure that
all selected ovals are in the boundary of the same tropical
component.
Take for~$\ind$ the
`monotonous'
function defined \via\ $\oind\oval=N+1-(\Dmax-\dpth\oval)$ if
$\oval\succcurlyeq\oval_i$ and
$\oind\oval=\Dmax-\dpth\oval$ if
$\oval\succcurlyeq\bar\oval_i$ for some~$i$, see
Figure~\ref{fig.index} (where the cases $N=7$ and $N=8$ are shown
schematically).
Due to Theorem~\ref{main} (see also Remark~\ref{rem.index}),
the pair $(C,\ind)$ is realized by a net of quadrics,
and for this net one has $b^0(V)=\omax{d}$,
see Lemma~\ref{imax=N}.
\qed

\midinsert
\centerline{\picture{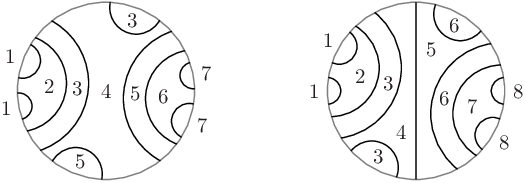}}
\figure
`Monotonous' index functions ($N=7$ and $N=8$)
\endfigure\label{fig.index}
\endinsert


\subsection{Proof of Theorem~\ref{th.B.2}}\label{proof.B.2}
In this section, we make an attempt to estimate the
Hilbert
number $\omax{d}$
introduced in
Definition~\ref{def.Hilbert}.

Let $C$ be a nonsingular real
plane algebraic curve of degree~$d$.
An oval of~$\CR$ is
said to be \emph{even} (\emph{odd}) if its depth is odd
(respectively, even). An oval is called \emph{hyperbolic} if it
has more than one immediate successor (in the partial order
defined in Section~\ref{s.topology}).
The following statement is known as
\emph{generalized Petrovsky inequality}.

\theorem[Theorem \rm(see~\cite{Arnold})]\label{th.Petrovsky}
Let~$C$ be a nonsingular real plane curve of even degree $d=2k$.
Then
$$
p-n^-\le\frac32k(k-1)+1,\quad
 n-p^-\le\frac32k(k-1).
$$
where $p$, $n$ are the numbers of even/odd
ovals of~$\CR$ and $p^-$, $n^-$ are the numbers of
even/odd hyperbolic ovals.
\qed
\endtheorem

\corollary\label{Petrovsky}
One has $\omax{d}\le\dfrac32k(k-1)+1$, where $k=[(d+1)/2]$.
\endcorollary

\proof
Let $d=2k$ be even, and let $C$ be a curve of degree~$d$ with
$m>1$ ovals of depth $k-1$. All submaximal ovals are
situated inside a nest
$\oval_1\prec\ldots\prec\oval_{k-2}$ of depth
$\Dmax-2=k-2$,
see Proposition~\ref{Bezout}. Assume that $k=2l$ is
even. Then the submaximal ovals are even, and one has
$p\ge m+l-1$, counting as well the even ovals
$\oval_1,\oval_3,\ldots\oval_{2l-1}$ in the nest.
On
the other hand, $n^-\le l-1$, as all odd ovals other than
$\oval_2,\oval_4,\ldots,\oval_{2l}$
are empty, hence not hyperbolic, see
Proposition~\ref{Bezout} again. Hence, the statement follows from
the first inequality in Theorem~\ref{th.Petrovsky}. The case of
$k$ odd is treated similarly, using the second inequality in
Theorem~\ref{th.Petrovsky}.

Let $d=2k-1$ be odd, and consider a real
curve~$C$ of degree~$d$ with a
nest $\oval_1\prec\ldots\prec\oval_{k-3}$ of depth $\Dmax-2=k-3$
and $m\ge2$ ovals $\oval',\oval'',\ldots$ of depth $k-2$. Pick a
pair of points $p'$ and~$p''$ inside~$\oval'$ and~$\oval''$,
respectively, and consider the line $L=(p_1p_2)$. From the B\'{e}zout
theorem, it follows that all points of intersection of~$L$ and~$C$
are one point on the one-sided component of~$\CR$ and a pair of
points on each of the ovals
$\oval_1,\ldots,\oval_{k-3},\oval',\oval''$. Furthermore, the pair
$\oval'$, $\oval''$ can be chosen so that all other innermost
ovals of~$\CR$ lie to one side of~$L$ in the interior
of~$\oval_{k-1}$ (which is divided by~$L$ into two components).
According to the Brusotti theorem~\cite{Brusotti},
the union $C+L$ can be perturbed to
form a nonsingular curve of degree~$2k$ with $m$ ovals of depth
$k-1$, see Figure~\ref{fig.2k-1} (where the curve and its
perturbation are shown schematically in grey and black, respectively).
Hence, the statement follows from the case of even degree
considered above.
\endproof

\midinsert
\centerline{\picture{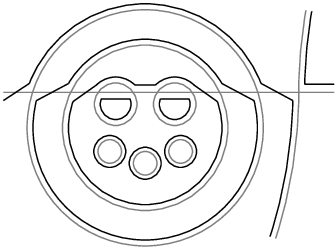}}
\figure
The perturbation of $C+L$ with a deep nest
\endfigure\label{fig.2k-1}
\endinsert

\theorem[Theorem \cite{Hilbertpaper}]\label{th.Hilbert}
For each
integer
$d\ge4$,
there is a nonsingular
curve of degree~$d$ in $\Rp2$
with
\Dashes
\dash
$4$ ovals of depth $1$ if $d = 4$,
\dash
$6$ ovals of depth $1$ if $d = 5$,
\dash
$k(k + 1) - 3$ ovals of
depth $[d/2]-1$ if $d = 2k$ is even,
\dash
$k(k + 2) - 3$ ovals
of depth $[d/2]-1$ if $d = 2k + 1$ is odd. \qed
\endDashes
\endtheorem

For the reader's convenience,
we give a brief outline of the original construction due to
Hilbert that
produces curves as in Theorem~\ref{th.Hilbert}.

For
even degrees,
one can use an inductive procedure
which produces a sequence of curves $C^{(2k)}$,
$\deg C^{(2k)}=2k$.
Let $E=\{p_E = 0\}$ be an ellipse in~$\Rp2$.
The curve~$C^{(2)}$ is defined
by a polynomial~$p^{(2)}$ of the form
$$
p^{(2)}= p_E + \varepsilon^{(2)}l^{(2)}_1 l^{(2)}_2
$$
where $\varepsilon^{(2)}>0$ is a real number,
$\ls|\varepsilon^{(2)}|\ll1$,
and $l^{(2)}_1$ and $l^{(2)}_2$ are real polynomials
of degree~$1$
such that the pair of lines $\{l^{(2)}_1 l^{(2)}_2 = 0\}$
intersects~$E$ at four distinct real points.
The intersection of the exterior of~$E$ and the interior of~$C^{(2)}$
is formed by two disks $D^{(2)}_1$ and $D^{(2)}_2$.
Inductively, we construct curves
$C^{(2k)}=\{
p^{(2k)}=0\}$
with the following
properties:
\roster
\item
$C^{(2k)}$ has an oval~$\oval^{(2k)}$ of depth~$k-1$
such that $\oval^{(2k)}$ intersects~$E$ at $4k$ distinct points,
the orders of the intersection points on $\oval^{(2k)}$ and $E$
coincide,
and the intersection
of the exterior of~$E$
and the exterior of~$\oval^{(2k)}$
consists of
a M\"{o}bius strip
and $2k-1$ discs $D^{(2k)}_1$, $\ldots$, $D^{(2k)}_{2k - 1}$
(shaded in Figure~\ref{fig.Hilbert}),
\item
one has
$$
p^{(2k)} = p^{(2k - 2)}p_E + \varepsilon^{(2k)}
l^{(2k)}_1 \ldots l^{(2k)}_{2k},
$$
where
$\varepsilon^{(2k)}$
is a real number,
$\ls|\varepsilon^{(2k)}|\ll1$,
and
$l^{(2k)}_1, \ldots, l^{(2k)}_{2k}$ are certain polynomials
of degree~$1$
such that the union of lines
$\{l^{(2k)}_1 \ldots l^{(2k)}_{2k}=0\}$ intersects~$E$ at $4k$
distinct real points, all points belonging to
$\partial D^{(2k - 2)}_1$.
\endroster
The sign of~$\varepsilon^{(2k)}$ is chosen
so
that
$\oval^{(2k - 2)}\cup E$ produces $4k - 4$ ovals of $C^{(2k)}$.

\midinsert
\centerline{\picture{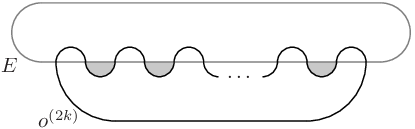}}
\figure
The oval $\oval^{(2k)}$ and the disks $D^{(2k)}_i$ (shaded)
\endfigure\label{fig.Hilbert}
\endinsert

The above properties imply that
each
curve $C^{(2k)}$ has the required number
of ovals of depth~$k - 1$
(and next curve~$C^{(2k + 2)}$ still satisfies \therosteritem1\,).

\midinsert
\centerline{\picture{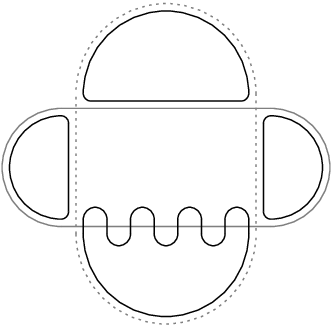}}
\figure
Hilbert's construction in degree $4$
\endfigure\label{fig-sextic}
\endinsert

The curves $C^{(2k + 1)}$ of odd degree $2k+1$ are constructed
similarly, starting from a curve~$C^{(3)}$ defined by a polynomial of
the form
$$
p^{(3)} = l p_E + \varepsilon^{(3)}l^{(3)}_1 l^{(3)}_2 l^{(3)}_3,
$$
where $\varepsilon^{(3)}>0$ is
a sufficiently small
real number,
$l$ is a polynomial of degree~$1$ defining a line
disjoint from~$E$,
and $l^{(3)}_1$, $l^{(3)}_2$, and $l^{(3)}_3$ are polynomials
of degree~$1$
such that the union of lines $\{l^{(3)}_1 l^{(3)}_2 l^{(3)}_3 = 0\}$
intersects~$E$ at six distinct real points.

\corollary\label{Hilbert}
One has $\omax{d}>\dfrac14(d-2)(d+4)-2$.
\qed
\endcorollary

\Remark
According to S.~Orevkov~\cite{Orevkov},
there are real algebraic
curves of degree~$d$ with
$$
\frac9{32}d^2+O(d)
$$
ovals of depth $[\frac{d}2]-1$, and it is expected that this
estimate is
still not sharp.
In the category of real pseudo-holomorphic curves, Orevkov achieved as many
as $\frac13d^2+O(d)$ ovals of submaximal depth.
\endRemark

\proof[Proof of Theorem~\ref{th.B.2}]
The statement of the theorem
follows from Theorem~\ref{th.estimate} and the bounds
on~$\omax{d}$ given by Corollaries~\ref{Petrovsky} and~\ref{Hilbert}.
\endproof

\section{Intersections of quadrics of dimension one}\label{S.curves}

In this section, we consider the case $N=r+2$, \ie,
one-dimensional complete intersections of quadrics.



\subsection{Proof of Theorem~\ref{th.B.-2}: the upper bound}
Let~$V$ be a regular complete intersection of $N - 1$
quadrics in $\CpN$.
Iterating the
adjunction formula,
one finds that
the genus $g(V)$
of the curve~$V$ satisfies
the relation
$$
2g(V)-2=2^{N-1}\bigl(2(N - 1)- (N + 1)\bigr)=2^{N - 1}(N - 3);
$$
hence, $g(V)=2^{N - 2}(N - 3) + 1$
and the Harnack inequality gives the upper bound
$\MAX(N-2,N) \leq 2^{N-2}(N-3)+2$.
\qed

\subsection{Proof of Theorem~\ref{th.B.-2}: the construction}
To prove the
lower bound
$\MAX(N-2,N)\ge 2^{N-2}(N-3)+2$,
for each integer $N\ge2$
we construct a
homogeneous quadratic
polynomial
$q\ix{N}\in \R[x_0, \ldots, x_N]$
and a
pair $(l\ix{N}_1, l\ix{N}_2)$
of linear forms $l\ix{N}_i \in \R[x_0, \ldots, x_N]$,
$i = 1$, $2$, with the following properties:
\roster
\item\local{curve.1}
the common zero set $V\ix{N}=\{q\ix2=\ldots=q\ix{N}=0\}\subset\CpN$
is a regular complete intersection;
\item\local{curve.2}
the real part $\VR\ix{N}$
has $2^{N - 2}(N - 3) + 2$ connected components;
\item\local{curve.3}
there is a distinguished component
$\oval\ix{N}\subset\VR\ix{N}$,
which has
two disjoint closed arcs $A\ix{N}_1$, $A\ix{N}_2$
such that the interior of $A\ix{N}_i$, $i = 1$, $2$,
contains all $2^{N-1}$
points
of intersection
of the hyperplane
$L\ix{N}_i=\{l\ix{N}_i=0\}$
with~$V\ix{N}$.
\endroster
Property~\loccit{curve.2} gives the desired lower bound.

The construction is by induction. Let
$$
l\ix2_1=x_2,\quad l\ix2_2=x_2-x_1,\quad\text{and}\quad
q\ix2=
l\ix2_1l\ix2_2
+(x_1-x_0)(x_1-2x_0).
$$
Assume that, for all integers $2\le k\le N$, polynomials
$q\ix{k}$, $l\ix{k}_1$, and $l\ix{k}_2$ satisfying
conditions \loccit{curve.1}--\loccit{curve.3} above are
constructed.
Let $\tilde l\ix{N}_2=l\ix{N}_2-\delta\ix{N}l\ix{N}_1$, where
$\delta\ix{N}>0$ is a real number so small that, for all
$t\in[0,\delta\ix{N}]$, the line
$\{l\ix{N}_2-t\ix{N}l\ix{N}_1=0\}$ intersects~$\oval\ix{N}$ at
$2^{N-1}$ distinct real points which all belong to the arc
$A\ix{N}_2$.
Put
$$
l\ix{N+1}_1=x_{N+1}\quad\text{and}\quad
l\ix{N+1}_2=x_{N+1}-l\ix{N}_1.
$$
The intersection of the cone
$\{q\ix2=\ldots=q\ix{N}=0\}\subset\Rp{N + 1}$
(over~$\VR\ix{N}$)
and the hyperplane
$L\ix{N+1}_2 = \{l\ix{N + 1}_2 = 0\}$
is a copy of $\VR\ix{N}$,
see Figure~\ref{fig-curve}.

\midinsert
\centerline{\picture{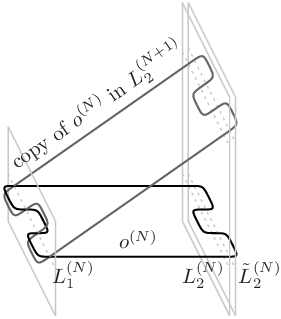}}
\figure
Construction of one-dimensional intersection of quadrics
\endfigure\label{fig-curve}
\endinsert

Put
$$
q\ix{N + 1} = l\ix{N + 1}_1 l\ix{N + 1}_2 +
\varepsilon\ix{N + 1}l\ix{N}_2\tilde{l}\ix{N}_2,
$$
where $\varepsilon\ix{N + 1}>0$ is a
sufficiently small
real number.
One can observe that,
on the hyperplane $\{l\ix{N}_1=0\}\subset\Rp{N+1}$, the polynomial
$q\ix{N + 1}$ has no
zeroes outside the
subspace $\{x_{N+1}=l\ix{N}_2=0\}$.

The new curve $V\ix{N + 1}$ is a regular complete intersection,
and its real part has $2^{N - 1}(N - 2) + 2$
connected components. Indeed, each component
$\oval\subset\VR\ix{N}$ other than $\oval\ix{N}$
gives rise to two components of $\VR\ix{N + 1}$, whereas
$\oval\ix{N}$ gives rise to $2^{N-1}$ components
of $\VR\ix{N + 1}$, each component being the perturbation of the
union $\arc_j\cup\arc_j'$, where $\arc_j\subset\oval\ix{N}$,
$j=0,\ldots,2^{N-1}-1$,
is the
arc bounded by two consecutive (in~$\oval\ix{N}$)
points of the intersection $L\ix{N}_1\cap\oval\ix{N}$ (and not
containing other intersection points)
and $\arc_j'$
is the copy of~$\arc_j$ in~$L\ix{N+1}_2$. All but one
of the
arcs~$\arc_j$ belong to $A\ix{N}_1$ and produce `small'
components; the arc $\arc_0$ bounded by the two outermost (from
the point of view of $A\ix{N}_1$) intersection points produces the
`long' component, which we take for $\oval\ix{N+1}$. Finally,
observe that the new component $\oval\ix{N+1}$ has two arcs
$A\ix{N+1}_i$, $i=1$, $2$,
satisfying condition~\loccit{curve.3} above: they
are the perturbations of the arc $A\ix{N}_2\subset L\ix{N+1}_1$
and its copy in $L\ix{N+1}_2$.
\qed


\section{Concluding remarks}\label{S.concluding}

In this section, we consider the few first special cases, $N=2$,
$3$, $4$, and~$5$, where, in fact, a complete deformation
classification can be given. We also briefly discuss the other
Betti numbers and the maximality of common zero sets of nets of
quadrics; however, we merely outline the directions of the further
investigation, leaving all details for a subsequent
paper.

\subsection{Empty intersections of quadrics}\label{s.class}
Consider a complete intersection~$V$ of $(r+1)$ real quadrics
in~$\CpN$, and assume that $\VR=\varnothing$. Choosing generators
$q_0,q_1,\ldots,q_r$ of the linear system, we obtain a map
$$
S^N=(\R^{N+1}\sminus0)/\R_+\to S^r=(\R^{r+1}\sminus0)/\R_+,\quad
u\mapsto(q_0(u),\ldots,q_r(u))/\R_+.
$$
Clearly, the
homotopy class of this map, which can be regarded as an element of
the group $\pi_N(S^r)$ modulo the antipodal involution,
is a deformation invariant
of the system.
Besides, the map is even (the images of~$u$ and $-u$ coincide);
hence, it
also induces certain maps $\RpN\to S^r$, $S^N\to\Rpr$, and
$\RpN\to\Rpr$, and
their homotopy classes are also deformation invariant.
Below, among other topics, we consider a few
special cases where
these classes do
distinguish empty regular intersections.

In general, the deformation classifications of linear
systems of quadrics, quadratic (rational) maps $S^N\to S^r$ (or
$\RpN\to\Rpr$), and spectral hypersurfaces (\eg, spectral curves,
even endowed with a theta characteristic) are different problems.
We will illustrate this by examples.

\subsection{Three conics}\label{s.Rp2}
We start with the case $r=N=2$, \ie, a net of conics in~$\Rp2$.
The spectral curve is a cubic $C\subset\Cp2$, and
the regularity condition implies that the common zero set must be
empty (even over~$\C$).
There are two deformation classes of
complete intersections of three conics; they can be distinguished
by the $\ZZ$-Kronecker invariant, \ie, the ${\bmod2}$ degree of the
associated map $\Rp2\to S^2$, see~\ref{s.class}.
If $\deg=0\bmod2$, the index function
takes all three values~$0$, $1$, and~$2$; otherwise, the
index function only takes the middle value~$1$.
(Alternatively, the Kronecker invariant counts
the parity of the number of real solutions of the system
$q_a=q_b=0$, $q_c>0$ in $\Rp2=\RpN$, where $a,b,c$ represent any
triple of
non-collinear points in $\Rp2=\Rpr$).

The classification of generic nets of conics can be obtained using
the results of~\cite{degt}. Note that, considering
generic quadratic maps $\Rp2\to\Rp2$ rather than
regular complete intersections,
there are four deformation
classes; they can be
distinguished by the topology of~$\CR$ and the spectral theta
characteristic.

\subsection{Spectral curves of degree~$4$}\label{s.Rp3}
Next special case is an
intersection of three quadrics in~$\Rp3$.
Here, a regular intersection~$\VR$ may consist of $0$, $2$, $4$, $6$, or $8$
real points, and the spectral curve is a quartic $C\subset\Rp2$.
Assuming $C$ nonsingular and computing the Euler characteristic (\eg, using
Theorem~\ref{th.ss} or the general formula for the Euler characteristic
found in~\cite{AG}) one can see that, if $V\ne\varnothing$, the
real part $\CR$ consists of $\frac12\Card\VR$ empty ovals. In this
case, $\Card\VR$ determines the net up to deformation.
If $\VR=\varnothing$, then either $\CR=\varnothing$ or $\CR$ is a
nest of depth two.
Such nets
form two deformation classes,
the homotopy class of the associated
quadratic map, see~\ref{s.class}, being either~$0$ or
$1\in\pi_3(S^2)/{\pm1}$.

\subsection{Canonical curves of genus~$5$ in~$\Rp4$}\label{s.Rp4}
Regular complete intersections of three quadrics in
$\Cp4$ are canonical curves of genus~$5$. Thus,
the set of projective classes of
such (real)
intersections is embedded
into the moduli space of
(real)
curves of genus~$5$. As is known, see,
\eg,~\cite{Tyurin}, the image of this embedding is the complement
of the strata
formed by the
hyperelliptic curves, trigonal curves, and curves
with a vanishing theta constant. Since each of
the three
strata
has
positive codimension,
the
known
classification of real forms of curves
of a given genus (applied to $g=5$) implies that the maximal number
of connected components
that
a regular complete intersection of three real quadrics in~$\Rp4$
may
have
equals~$6=\omax5$.


\subsection{$K3$-surfaces of degree~$8$ in~$\Rp5$}\label{s.Rp5}
A regular complete intersection of three quadrics in the
projective space of dimension~$5$ is a $K3$-surface with a
(primitive) polarization of degree~$8$. Thus, as in the previous
case, the set of projective classes of intersections is embedded
into the moduli space of $K3$-surfaces with a polarization of
degree $8$,
the complement
consisting of a few strata of positive codimension (for details,
see~\cite{SD}). In particular, any generic $K3$-surface with a
polarization of degree~$8$ is indeed a complete intersection of
three quadrics. The deformation classification of $K3$-surfaces
can be obtained using the results of V.~V.~Nikuin~\cite{Nikulin}.
The case of maximal real $K3$-surfaces is
particulary
simple: there are three deformation classes, distinguished by the
topology of the real part, which can be $S_{10}\sqcup S$,
$S_6\sqcup 5S$, or $S_2\sqcup 9S$. In particular, the maximal
number of connected components of a complete intersection of three
real quadrics in~$\Rp5$ equals $10=\omax{6}+1$. (There is another
shape with $10$ connected components, the $K3$-surface with the
real part $S_1\sqcup 9S$. However, one can easily show that a
$K3$-surface of degree~$8$ cannot have $10$ spheres.)

\subsection{Other Betti numbers}\label{s.Betti}
The techniques of this paper
can be used to estimate the other Betti numbers as well.
For $0\le i<\frac12(N-3)$, it would give a bound of the form
$$
B^i_2(N)-\Hilb_{i+1}(N+1)=O(1),
$$
where $B^i_2(N)$ is
the maximal $i$-th Betti number
of a
regular complete intersection of three real quadrics in $\RpN$
and $\Hilb_{i+1}(N+1)$ is the maximal number of ovals of depth
$\ge(\Dmax-i-1)=[\frac12(N-1)]-i$ that a nonsingular real plane
curve of degree $d=N+1$ may have. The possible discrepancy
is due to
a couple of
unknown
differentials in the spectral sequence and the inclusion
homomorphism $H^{N-3-i}(\RpN)\to H^{N-3-i}(\VR)$.

\subsection{The examples are asymptotically
maximal}\label{s.maximality}
Recall that, given a real algebraic variety~$X$,
the Smith inequality
states that
$$
\dim H_*(X_\R)\le\dim H_*(X).
\eqtag\label{eq.Smith}
$$
(As usual, all homology groups are with the $\ZZ$ coefficients.)
If the equality holds, $X$ is said to be \emph{maximal},
or an \emph{$M$-variety}.
In particular, if $X$ is a nonsingular plane curve of degree
$N+1$, the Smith inequality~\eqref{eq.Smith}
implies that the number of connected
components of $X_\R$ does not exceed $g+1=\frac12N(N-1)+1$.
The Hilbert curves
used
in Section~\ref{topology.zerolocus}
to construct nets with a big number of
connected components are
known to be maximal.
Using the spectral sequence of Theorem~\ref{th.ss}, one can easily
see that,
under the choice of the index
function made in the proof of Theorem~\ref{th.estimate},
the dimension $\dim H_*(\Lm)$
equals
$2(g+1)+2$ if $N$ is odd or
$2(g+1)+1$ if $N$ is even.

On the other hand, for the common zero set~$V$ (as for any
projective variety) there is a certain constant~$l$
such that the inclusion
homomorphism $H_i(\VR)\to H_i(\RpN)$ is nontrivial for all
$i\le l$ and trivial for all $i>l$
(see~\cite{Kh};
in our case $1\le l<\frac12N$).
Hence, by the Poincar\'{e}--Lefschetz duality and
Lemma~\ref{L+=complement}, one has
$$
\dim H_*(\Lm)=\dim H_*(\RpN,\VR)=\dim H_*(\VR)+N-2l-1.
$$
Finally,
one can easily find $\dim H_*(V)$:
it equals $4(k^2-1)$ if $N=2k$ is even and $4(k^2-k)$ if $N=2k-1$
is odd.
Now,
combining the above computation,
one observes that the intersections of quadrics constructed from
the Hilbert curves using monotonous index functions are
asymptotically maximal in the sense that
$$
\dim H_*(\VR)=\dim H_*(V)+O(N)=N^2+O(N).
$$

The latter identity
shows that the upper bound for $\MAX(2,N)$ provided by
the Smith inequality is
too
rough: this bound is of the
form
$\frac12N^2+O(N)$,
whereas, as is shown in this paper,
$\MAX(2,N)$ does not exceed $\frac38N^2+O(N)$.
When the intersection is of even dimension,
one can improve the leading coefficient
in the bound
by combining the Smith inequality
and
the generalized Comessatti
inequality;
however, the resulting estimate is still
too far from the sharp bound.

\subsection{The examples are not maximal}\label{s.not.maximal}
Another interesting consequence
of the computation of the previous section is the fact that,
starting from $N=6$, the complete intersections of quadrics
maximizing the number of components are never truly maximal in the
sense of the Smith inequality~\eqref{eq.Smith}: one has
$$
N+O(1)\le\dim H_*(V)-\dim H_*(\VR)\le2N+O(1).
$$
Using the spectral sequence of Theorem~\ref{th.ss}, one can
easily show
that a maximal complete intersection~$V$ of three real quadrics
in~$\RpN$ must have
index function taking values between $\frac12(N-1)$ and
$\frac12(N+3)$
(\cf.~\eqref{imax}--\eqref{qmax});
the real part~$\VR$ has large Betti numbers in two
or three middle dimensions, (most) other Betti numbers being equal
to~$1$.

Apparently, it is the Harnack $M$-curves that are suitable for
obtaining nets with maximal common zero locus. However, at present
we do not know much about the differentials in the spectral
sequence or the constant~$l$ introduced in the previous section.
It may happen that these data are controlled by
an extra
flexibility
in the choice of the real $\Spin$-structure on the spectral
curve:
in addition to the semi-orientation, one can also choose the
values on the components of~$\CR$.

\refstyle{C}

\enddocument